\documentclass{article}
\usepackage{amssymb,amsmath,
theorem,euscript}

\newcounter{sec}

\newcounter{punct}[sec]

\def\punct{\refstepcounter{punct}{\arabic{sec}.\arabic{punct}.  }}

\def\COUNTERS{\addtocounter{sec}{1}
              \setcounter{punct}{0}
          \setcounter{equation}{0}
          \setcounter{theorem}{0}
          \setcounter{problem}{0}
            \setcounter{Apunct}{0}
          }

\def\COUNTERS{\addtocounter{sec}{1}
              \setcounter{punct}{0}
          \setcounter{equation}{0}
          \setcounter{theorem}{0}
          }

\begin{document}

%%%1. ClASSICAL GROUPS
\def\SL{\mathrm {SL}}
\def\SU{\mathrm {SU}}
\def\GL{\mathrm  {GL}}
\def\U{\mathrm  U}
\def\OO{\mathrm  O}
\def\Sp{\mathrm  {Sp}}
\def\SO{\mathrm  {SO}}
\def\SOS{\mathrm {SO}^*}

\def\PGL{\mathrm  {PGL}}
\def\PU{\mathrm {PU}}

\def\Gr{\mathrm{Gr}}

\def\Fl{\mathrm{Fl}}

\def\OSp{\mathrm {OSp}}

\def\Mat{\mathrm{Mat}}

\def\Pfaff{\mathrm {Pfaff}}

\def\B{\mathbf B}

\def\phi{\varphi}
\def\epsilon{\varepsilon}
\def\kappa{\varkappa}

\def\le{\leqslant}
\def\ge{\geqslant}

\renewcommand{\Re}{\mathop{\rm Re}\nolimits}

\renewcommand{\Im}{\mathop{\rm Im}\nolimits}

\def\pia{\pi_\downarrow}

\newcommand{\im}{\mathop{\rm im}\nolimits}
\newcommand{\indef}{\mathop{\rm indef}\nolimits}
\newcommand{\dom}{\mathop{\rm dom}\nolimits}
\newcommand{\codim}{\mathop{\rm codim}\nolimits}

\def\cA{\mathcal A}
\def\cB{\mathcal B}
\def\cC{\mathcal C}
\def\cD{\mathcal D}
\def\cE{\mathcal E}
\def\cF{\mathcal F}
\def\cG{\mathcal G}
\def\cH{\mathcal H}
\def\cJ{\mathcal J}
\def\cI{\mathcal I}
\def\cK{\mathcal K}
\def\cL{\mathcal L}
\def\cM{\mathcal M}
\def\cN{\mathcal N}
\def\cO{\mathcal O}
\def\cP{\mathcal P}
\def\cQ{\mathcal Q}
\def\cR{\mathcal R}
\def\cS{\mathcal S}
\def\cT{\mathcal T}
\def\cU{\mathcal U}
\def\cV{\mathcal V}
\def\cW{\mathcal W}
\def\cX{\mathcal X}
\def\cY{\mathcal Y}
\def\cZ{\mathcal Z}

\def\frA{\mathfrak A}
\def\frB{\mathfrak B}
\def\frC{\mathfrak C}
\def\frD{\mathfrak D}
\def\frE{\mathfrak E}
\def\frF{\mathfrak F}
\def\frG{\mathfrak G}
\def\frH{\mathfrak H}
\def\frJ{\mathfrak J}
\def\frK{\mathfrak K}
\def\frL{\mathfrak L}
\def\frM{\mathfrak M}
\def\frN{\mathfrak N}
\def\frO{\mathfrak O}
\def\frP{\mathfrak P}
\def\frQ{\mathfrak Q}
\def\frR{\mathfrak R}
\def\frS{\mathfrak S}
\def\frT{\mathfrak T}
\def\frU{\mathfrak U}
\def\frV{\mathfrak V}
\def\frW{\mathfrak W}
\def\frX{\mathfrak X}
\def\frY{\mathfrak Y}
\def\frZ{\mathfrak Z}

\def\fra{\mathfrak a}
\def\frb{\mathfrak b}
\def\frc{\mathfrak c}
\def\frd{\mathfrak d}
\def\fre{\mathfrak e}
\def\frf{\mathfrak f}
\def\frg{\mathfrak g}
\def\frh{\mathfrak h}
\def\fri{\mathfrak i}
\def\frj{\mathfrak j}
\def\frk{\mathfrak k}
\def\frl{\mathfrak l}
\def\frm{\mathfrak m}
\def\frn{\mathfrak n}
\def\fro{\mathfrak o}
\def\frp{\mathfrak p}
\def\frq{\mathfrak q}
\def\frr{\mathfrak r}
\def\frs{\mathfrak s}
\def\frt{\mathfrak t}
\def\fru{\mathfrak u}
\def\frv{\mathfrak v}
\def\frw{\mathfrak w}
\def\frx{\mathfrak x}
\def\fry{\mathfrak y}
\def\frz{\mathfrak z}

\def\fros{\mathfrak{s}}

\def\bfa{\mathbf a}
\def\bfb{\mathbf b}
\def\bfc{\mathbf c}
\def\bfd{\mathbf d}
\def\bfe{\mathbf e}
\def\bff{\mathbf f}
\def\bfg{\mathbf g}
\def\bfh{\mathbf h}
\def\bfi{\mathbf i}
\def\bfj{\mathbf j}
\def\bfk{\mathbf k}
\def\bfl{\mathbf l}
\def\bfm{\mathbf m}
\def\bfn{\mathbf n}
\def\bfo{\mathbf o}
\def\bfp{\mathbf q}
\def\bfr{\mathbf r}
\def\bfs{\mathbf s}
\def\bft{\mathbf t}
\def\bfu{\mathbf u}
\def\bfv{\mathbf v}
\def\bfw{\mathbf w}
\def\bfx{\mathbf x}
\def\bfy{\mathbf y}
\def\bfz{\mathbf z}

\def\bfA{\mathbf A}
\def\bfB{\mathbf B}
\def\bfC{\mathbf C}
\def\bfD{\mathbf D}
\def\bfE{\mathbf E}
\def\bfF{\mathbf F}
\def\bfG{\mathbf G}
\def\bfH{\mathbf H}
\def\bfI{\mathbf I}
\def\bfJ{\mathbf J}
\def\bfK{\mathbf K}
\def\bfL{\mathbf L}
\def\bfM{\mathbf M}
\def\bfN{\mathbf N}
\def\bfO{\mathbf O}
\def\bfP{\mathbf P}
\def\bfQ{\mathbf Q}
\def\bfR{\mathbf R}
\def\bfS{\mathbf S}
\def\bfT{\mathbf T}
\def\bfU{\mathbf U}
\def\bfV{\mathbf V}
\def\bfW{\mathbf W}
\def\bfX{\mathbf X}
\def\bfY{\mathbf Y}
\def\bfZ{\mathbf Z}

\def\R {{\mathbb R }}
 \def\C {{\mathbb C }}
  \def\Z{{\mathbb Z}}
  \def\H{{\mathbb H}}
\def\K{{\mathbb K}}
\def\N{{\mathbb N}}
\def\Q{{\mathbb Q}}
\def\A{{\mathbb A}}

\def\T{\mathbb T}

\def\bbA{\mathbb A}
\def\bbB{\mathbb B}
\def\bbD{\mathbb D}
\def\bbE{\mathbb E}
\def\bbF{\mathbb F}
\def\bbG{\mathbb G}
\def\bbI{\mathbb I}
\def\bbJ{\mathbb J}
\def\bbL{\mathbb L}
\def\bbM{\mathbb M}
\def\bbN{\mathbb N}
\def\bbO{\mathbb O}
\def\bbP{\mathbb P}
\def\bbQ{\mathbb Q}
\def\bbS{\mathbb S}
\def\bbT{\mathbb T}
\def\bbU{\mathbb U}
\def\bbV{\mathbb V}
\def\bbW{\mathbb W}
\def\bbX{\mathbb X}
\def\bbY{\mathbb Y}

 \def\ov{\overline}
\def\wt{\widetilde}
\def\wh{\widehat}

\def\P{\mathbb P}

\def\arr{\rightrightarrows}

\def\SS{\smallskip}

\def\ev{{\mathrm{even}}}
\def\od{{\mathrm{odd}}}

\def\q{\quad}

\def\F{\mathbf F}

\def\b{\mathbf b}

\def\RA{\Longrightarrow}

\begin{center}

\bf\Large On beta-function of tube of light
cone

\sc\large

\bigskip

Yuri A. Neretin%
\footnote{Supported by the grant FWF, project P19064,
 Russian Federal Agency for Nuclear Energy,
  grant NWO.047.017.015, and grant JSPS-RFBR-07.01.91209 }

\end{center}

{\small We construct $B$-function of the Hermitian
symmetric space $\OO(n,2)/\OO(n)\times \OO(2)$
or equivalently of the tube
$(\Re z_0)^2> (\Re z_1)^2+\dots+ (\Re z_n)^2$
in $\C^{n+1}$}.

\bigskip

\section{Formulation of the result}

\COUNTERS

{\bf\punct Preliminary references.}
The beta-function of symmetric cones
$$
\GL(n,\R)/\OO(n), \,\, \GL(n,\C)/\U(n),\,\,
\GL(n,\H)/\Sp(n)
$$
 was constructed by Gindikin in \cite{Gin},
see also \cite{FK}.
For the remaining series of classical symmetric spaces
the beta-function was obtained in
\cite{Ner-beta}. The  subseries $\OO(n,2)/\OO(n)\times\OO(2)$
has two beta-functions, the first one is a special case
of the beta-function of
 $\OO(p,q)/\OO(p)\times\OO(q)$.
The second beta-function is discussed here, it is related
 to the Hermitian structure of these spaces.%
\footnote{This note is some kind of an addendum to my paper 
\cite{Ner-beta} on matrix beta-functions;
originally it was an omitted  section 
of the last work. Since the Plancherel formula for Berezin
representations (which was 
my actual purpose) for $\OO(n,2)/ \OO(n)\times\OO(2)$
 was known due Berezin himself \cite{Ber} (proof was published by
Unterberger and Upmeier \cite{UU}) I omitted 6-parameter
 integral (\ref{formula}). 
   However, after a discussion with prof. G.Roos I understood that this
integral can be interesting by itself since this series of symmetric spaces
is familiar to various non-representation-theoretic people.}

%%%%%%%%%%%%%%%%%%%%%%%%%%%%%%%%%%%%%%%%%%%%%%

\SS

 {\bf \punct The tube of light cone.}
Consider the space $\C^{n+1}$ with coordinates
$z_0$, $z_1$, \dots, $z_{n}$. By
$\frT_n$ we denote the tube
$$
(\Re z_0)^2 > (\Re z_1)^2 + (\Re z_2)^2+\dots+
 (\Re z_{n})^2,\qquad \Re z_0>0
$$
The space $\frT_n$ is homogeneous
with respect to the pseudo-orthogonal group
$\OO(n+1,2)$ (apparently, this was discovered
by E.Cartan,  \cite{Car}, for a discussion, see, for instance,
\cite{Hua}, \cite{Pya}, \cite {FK}), namely
$$
\frT_n\simeq\OO(n+1,2)/\OO(n+1)\times\OO(2)
$$
This group acts on $\frT $ by quadratic-fractional
transformations and is generated by the maps
of the following 3 types

\SS

a) Shifts $z\mapsto z+ia$, where $a\in\R^{n+1}$.

\SS

b) Linear transformations $z\mapsto zg$, where
$g\in\SO_0(n,1)$, i.e., $g$ is a real matrix
preserving  the quadratic form $x_0^2-x_1^2-\dots - x_n^2$.

\SS

c) Quasiinversion $z\mapsto z/(z,z)$, where
$(z,z):=z_0^2-z_1^2-\dots-z_n^2$.

\SS

This space has numerous names, in particular,
{\it Cartan domain of 4-th type},
{\it future tube}, {\it Lie spheres}, {\it Lie balls}%
\footnote{The last two terms arose as a
result of non-correct successive translations
German--Chines--English--Russian--English. 

Initially,
a {\it 'Lie sphere'} was  an oriented subsphere in $S^n=\R^n\cup\infty$
or a point. {\it Lie spheres} was the space of Lie spheres;
It is a homogeneous space  $\OO(n+2)/\OO(n+2)\times \OO(2)$.
{\it Lie sphere geometry} , see for instance \cite{Cec},
was a geometry of this space. 
Hua Loo Keng extended the 
term 'Lie spheres' to the dual symmetric space
$\OO(n+2,2)/\OO(n+2)\times \OO(2)$.

In Russian edition of Hua's book
\cite{Hua}
 the space
$\OO(n,2)/\OO(n)\times \OO(2)$ became
a {\it Lie sphere} ('sfera Li'; in fact  there are no traces of
the original meaning of 'Lie spheres' in Hua's book).
However, the term 'sphere' for an 
open domain is too peculiar and our space turned
into a {\it 'Lie ball'} in the English edition.}.

\SS

We prefer another realization of the same space.
Namely, consider $\C^{n+1}$ with coordinates
$u_1$, $u_2$, $z_1$, \dots, $z_{n-1}$,
\begin{align*}
u_1=v_1+i w_1;
\\
 u_2=v_2+iw_2;
 \\
 z_j:=x_j+y_j
\end{align*}

We define our tube $\frT_n$ by the inequalities
$$v_1v_2-\sum x_j^2> 0,\qquad v_1>0  $$

%%%%%%%%%%%%%%%%%%%%%%%%%%%%%%%%%%%%%%%%%%%%%%%

\SS

{\bf\punct Formulation of the result.}
For a nonzero complex $a$ and complex
$\lambda$, $\mu$ denote
$$
a^{\{\lambda|\mu\}}:=a^\lambda \ov a^\mu
$$
We also denote
$$
du_1:=dv_1\,dw_1,\quad dz_j:=dx_j\,dy_j,
\quad du:=du_1\,du_2,
\quad
dz:=dx_1\,dy_1\dots dx_{n-1}\,dy_{n-1}
$$
In particular $du\,dz$ denotes an integration with respect
to the Lebesgue measure on $\frT_n$.

Our purpose is the following formula
\begin{multline}
\int_{\frT_n}
\frac{
v_1^{\lambda_1-\lambda_2}(v_1v_2-\sum x_j^2)^{\lambda_2-n-1}
du\,dz
}
{(1+u_1)^{\{\sigma_1-\sigma_2|\tau_1-\tau_2\}}
\,
\bigl((1+u_1)(1+u_2)-\sum z_j^2\bigr)^{\{\sigma_2|\tau_2\}}
}
=\\=
 2^{1-\sigma_1-\tau_1+n}
\pi
\frac{
    \Gamma\left(\lambda_1-\tfrac12 (n+1)\right)
    \Gamma\left(\sigma_1+\tau_1-\lambda_1-\tfrac 12(n-1)\right)
    }
    {
    \Gamma\left(\sigma_1-\tfrac12 (n-1)\right)
     \Gamma\left(\tau_1-\tfrac12 (n-1)\right)
       }
\times \\ \times
2^{2-\sigma_2-\tau_2+n}
\pi^n
\frac{\Gamma(\lambda_2-n)\Gamma(\sigma_2+\tau_2-\lambda_2)}
{\Gamma(\sigma_2)\Gamma(\tau_2)}
\label{formula}
\end{multline}
where $\lambda_1$, $\lambda_2$,
$\sigma_1$, $\sigma_2$, $\tau_1$, $\tau_2\in\C$.

\SS

{\sc Remark.} Let us explain the meaning of complex powers.
The base numbers $v_1$, $v_1 v_2-\sum x_j^2$ of the numerator
 are positive reals. Next, the point $u_1=1$, $u_2=1$,
$z=0$ is contained in $\frT_n$. The denominator
is well defined at this point. Since the domain $\frT_n$ is simply connected,
the corresponding branches of the power functions are well defined.

\SS

%%%%%%%%%%%%%%%%%%%%%%%%%%%%%%%%%%%%%%%%

{\bf\punct Meaning of factors.}
a) Functions
$v_1^{\lambda_1-\lambda_2}(v_1v_2-\sum x_j^2)^{\lambda_2-n-1}$
are precisely the eigenfunctions of the parabolic
subgroup in $\OO(n+1,2)$.
Also $(v_1v_2-\sum x_j^2)^{-n-1}$
is the density of the $\OO(n+1,2)$-invariant
measure on $\frT_n$.

\SS

b) $\bigl((1+u_1)(1+u_2)-\sum z_j^2\bigr)$
is the standard term that is present
in formulae for the Cauchy kernel, the Bargman kernel,
see \cite{Hua}
and (more generally) for the Berezin kernels
on $\frT_n$.

\SS

%%%%%%%%%%%%%%%%%%%%%%%%%%%%%%%%%%%%%%%%%%%%%%

{\bf\punct  Comments. Special cases.}
a) If $\lambda_1=\lambda_2=0$, $\sigma_1=\tau_1=0$,
$\sigma_2=\tau_2$,
then we get one of Hua integrals, \cite{Hua}

\SS

b) The Plancherel formula
for Berezin kernels for the spaces
$\OO(n+1,2)/\OO(n+1)\times\OO(2)$
is reduced to our integral with
$\sigma_1=\tau_1=0$, $\sigma_2=\tau_2$ (see \cite{UU}). 
Apparently, Berezin himself (he perished in an accident in 1980)
derived this formula in some another way (however, his proof is
unknown; see also \cite{UU}).

\SS

c) For $n=3$, $4$, $6$ there are the following exceptional
isomorphisms of homogeneous spaces
\begin{align*}
\OO(3,2)/\OO(3)\times\OO(2)
&=
\Sp(6,\R)/\U(2);
\\
\OO(4,2)/\OO(4)\times\OO(2)
&=
\U(2,2)/\U(2)\times \U(2)
\\
\OO(6,2)/\OO(6)\times\OO(2)
&=
\SOS(8)/\U(4)
\end{align*}
In these cases,
our integrals coincide with matrix
beta-integrals obtained in \cite{Ner-beta}.

%%%%%%%%%%%%%%%%%%%%%%%%%%%%%%%%%%%%%%%%%%%%%%%%%%%%%%%%%%%%%%

\section{Calculations}

\COUNTERS

{\bf\punct Change of variables.}
Firstly, we transform our integral
to
$$
\int\limits_{v_1>0,\,\, w_1\in\R}
\frac{v_1^{\lambda_1-n-1}}{(1+u_1)^{\{\sigma_1|\tau_1\}}
}
\cdot
\Biggl\{\,\,
\int\limits_{v_2-\frac 1{v_1}\sum x_j^2>0}
\frac
{
\left(
v_2-\frac 1{v_1}\sum x_j^2
\right)^{\lambda_2-n-1}
du_2\,dz
}
{\left(1+u_2-\frac 1{1+u_1}\sum z_j^2\right)
^{\{\sigma_2|\tau_2\}}
}
\Biggr\} du_1
$$
Next, we change the variable
$v_2$ by $r$,
$$
r:=v_2-\frac 1{v_1}\sum x_j^2
$$
(the Jacobian of this substitution $=1$).
The interior integral now is reduced to
\begin{equation}
\int\limits_{r>0,\, w_1\in\R,\, x\in\R^{n-1},\,y\in\R^{n-1}}
\frac{r^{\lambda_2-n-1}
dr\, dw_2\,dx\,dy
}
{\left(
1+r+iw_2+\frac 1{v_1}\sum x_j^2-
\frac 1{1+u_1}\sum z_j^2\right)
^{\{\sigma_2|\tau_2\}}
}
\label{eq:after-first}
\end{equation}

{\bf\punct Next change of variables.}
Now we wish to decompose the expression
$$
H:=
1+r+iw_1+\frac 1{v_1}\sum x_j^2-
\frac 1{1+u_1}\sum z_j^2
$$
 as a sum
of imaginary and real parts.
For this purpose, we write
$$
\Re \frac 1{1+u_1}\sum z_j^2
=
\Re \frac{\sum (x_j^2- y_j^2+2 i x_j y_j)}
{1+v_1+iw_1}
=\sum_j\frac{(x_j^2-y_j^2)(1+v_1)+2 x_j y_jw_1}
{(1+v_1)^2+w_1^2}
$$
Therefore,
$$
\Re H=1+r+\sum_j
\begin{pmatrix} x_j& y_j\end{pmatrix}
S
\begin{pmatrix} x_j\\ y_j\end{pmatrix}
$$
where
$$
S:=
\begin{pmatrix}
\frac 1{v_1}-\frac{1+v_1}{(1+v_1)^2+w_1^2}
               &
    \frac{-w_1}{(1+v_1)^2+w_1^2}\\
    \frac{-w_1}{(1+v_1)^2+w_1^2}&
    \frac{1+v_1}{(1+v_1)^2+w_1^2}
    \end{pmatrix}
$$
Note that,
$$
\det S=\frac 1{v_1\bigl((1+v_1)^2+w_1^2\bigr)}
=v_1^{-1}(1+v_1+iw_1)^{-\{1|1\}}
=
v_1^{-1}(1+u_1)^{-\{1|1\}}
$$
Thus, $\det S>0$, the diagonal elements of $S$
also are positive. Hence $S$ is positive definite;
therefore $S^{1/2}$ is well defined.

\SS

Our next change of variables is
$$
\begin{pmatrix} x_j& y_j\end{pmatrix}
=\begin{pmatrix} p_j& q_j\end{pmatrix}S^{-1/2}
$$
Its Jacobian is
$$
v_1^{(n-1)/2} (1+u_1)^{(n-1)/2\,\cdot\,\{1|1\}}
$$
Finally, the 'interior integral' (\ref{eq:after-first})
comes to the form
\begin{multline*}
v_1^{(n-1)/2} (1+u_1)^{(n-1)/2\,\cdot\,\{1|1\}}
\times\\
\times
\int_{r>0, w_2\in\R, p\in \R^{n-1},\,q\in \R^{n-1}}
\frac{r^{\lambda_2-n-1}\,dr \,dw_2\,dp\,dq}
{\left(1+r+\sum(p_j^2+q_j^2)
+ iw_2 +i Q(p,q, v_1, w_1)\right)
^{\{\sigma_2|\tau_2\}}
}
\end{multline*}
where
$Q(p,q, v_1, w_1)$ is a real expression
(its explicit form is unessential  for us).

\SS

%%%%%%%%%%%%%%%%%%%%%%%%%%%%%%%%%%%%%%%%%%%%%%%

{\bf\punct Separation of variables.}
Now we change the variable
$w_2$ to
$$h:=w_2+Q(\cdot)$$
 The Jacobian is 1,
and we reduce our initial integral to the product
\begin{multline}
I\cdot J:=
\int_{v_1>0,w_1\in\R}
\frac{v_1^{\lambda_1-n/2-3/2}dv_1\,dw_1}
{(1+v_1+i w_1)^{\{\sigma_1-(n-1)/2|\tau_1-(n-1)/2\}}
}
\times
\\ \times
\int_{r>0, \, h\in\R, \,p\in\R^{n-1},\,q\in\R^{n-1}}
\frac{
r^{\lambda_2-n-1} dr\,d h \,dp\, dq
}
{\left(1+r+ \sum(p_j^2+q_j^2)+ i h\right)
^{\{\sigma_2|\tau_2\}}
}
\label{eq:separated}
\end{multline}
where $I$ denotes the first integral factor,
and $J$ the second one.

\SS

%%%%%%%%%%%%%%%%%%%%%%%%%%%%%%%%%%%%%%%%%%

{\bf\punct An auxiliary integral.}
First, we derive the identity
\begin{equation}
\int_{x>0,y\in\R}\frac{x^{\alpha-1} dx\,dy}
{(1+x+iy)^{\{\beta|\gamma\}}}
= 2^{2-\beta-\gamma}\pi
\frac{\Gamma(\alpha)\Gamma(\beta+\gamma-\alpha-1)}
         {\Gamma(\beta)\Gamma(\gamma)}
\label{eq:2-beta}
\end{equation}

We represent the left-hand side as
$$
\int_0^\infty dv\cdot v^{\alpha-1}
\int_{-\infty}^\infty
 \frac {dw} {(1+v+iw)^{\{\beta|\gamma\}}}
$$

The interior integral is the Cauchy beta-integral
(see \cite{PBM}, 2.2.6.31) and we get
$$
2\pi\cdot 2^{1-\beta-\gamma}
\frac{\Gamma(\beta+\gamma-1)}
{\Gamma(\beta)\Gamma(\gamma)}
\int_0^\infty \frac{v^{\alpha-1}dv}{(1+v)^{\beta+\gamma-1}}
$$

The last integral is a rephrasing of the definition
of the beta-function (see \cite{PBM}, 2.2.4.29).
We get
$$
2\pi\cdot 2^{1-\beta-\gamma}
\frac{\Gamma(\beta+\gamma-1)}
{\Gamma(\beta)\Gamma(\gamma)}
\cdot
\frac{\Gamma(\alpha)\Gamma(\beta+\gamma-\alpha-1)}
{\Gamma(\beta+\gamma-1)}
$$

{\bf\punct The first factor.} By (\ref{eq:2-beta}),
\begin{equation}
I:= 2^{1-\sigma_1-\tau_1+n}
\pi
\frac{
    \Gamma\left(\lambda_1-\tfrac12 (n+1)\right)
    \Gamma\left(\sigma_1+\tau_1-\lambda_1-\tfrac 12(n-1)\right)
    }
    {
    \Gamma\left(\sigma_1-\tfrac12 (n-1)\right)
     \Gamma\left(\tau_1-\tfrac12 (n-1)\right)
       }
       \label{eq:I}
\end{equation}

%%%%%%%%%%%%%%%%%%%%%%%%%%%%%%%%%%%%%%%%%%%%

\SS

{\bf \punct The second factor.}
Now we evaluate the factor $J$
in (\ref{eq:separated}). First, we pass
to the spherical coordinates
$$
R^2=\sum(p_j^2+q_j^2)
$$
in $\R^{2n-2}$. We get
(see \cite{PBM}, 3.3.2.1)
$$
\frac{2\pi^{n-1}}
{\Gamma(n-1)}
\int_{r>0,\, R>0,\,h\in\R}
\frac{r^{\lambda_2-n-1}R^{2n-3}dr\,dR\,dh}
{(1+r+R^2+ih)^{\{\sigma_2|\tau_2\}}}
$$
Next, we substitute
$\rho:=R^2$,
$$
J=\frac{\pi^{n-1}}
{\Gamma(n-1)}
\int_{r>0,\, R>0,\,h\in\R}
\frac{r^{\lambda_2-n-1} \rho^{n-2}\,dr\,d\rho\,dh}
{(1+r+\rho+ih)^{\{\sigma_2|\tau_2\}}}
$$
Next, we pass from the variables $(r,\rho)$
to
$$
(x,r):=(r+\rho,r)
$$
i.e.,
$$
J=\frac{\pi^{n-1}}
{\Gamma(n-1)}
\int_{x>0,\, 0<r<x,\,h\in\R}
\frac{r^{\lambda_2-n-1} (x-r)^{n-2}\,dr\,dx\,dh}
{(1+x+ih)^{\{\sigma_2|\tau_2\}}}
$$
Now we integrate in
$r$ using the standard definition of the beta-function,
$$
J=\frac{\pi^{n-1}}
{\Gamma(n-1)}
\mathrm B
(\lambda-n,n-1)
\int_{x>0,\,h\in\R}
\frac{x^{\lambda_2-2}\,dx\,dh}
{(1+x+ih)^{\{\sigma_2|\tau_2\}}}
$$
The last integral is of the form (\ref{eq:2-beta}).
Finally,
$$
J=2^{2-\sigma_2-\tau_2-n}
\pi^n
\frac{\Gamma(\lambda_2-n)\Gamma(\sigma_2+\tau_2-\lambda_2)}
{\Gamma(\sigma_2)\Gamma(\tau_2)}
$$

{\tt Math.Dept., University of Vienna,

 Nordbergstrasse, 15,
Vienna, Austria

\&

Institute for Theoretical and Experimental Physics,

Bolshaya Cheremushkinskaya, 25, Moscow 117259,
Russia

e-mail: neretin(at) mccme.ru

URL:www.mat.univie.ac.at/$\sim$neretin
}

\end{document}